\documentclass[11pt]{amsart}
\usepackage[english]{babel}
\usepackage{graphicx}
\usepackage{framed}
\usepackage[normalem]{ulem}
\usepackage{amsmath}
\usepackage{amsthm}
\usepackage{amssymb}
\usepackage{amsfonts}
\usepackage{enumerate}
\usepackage{hyperref}
\usepackage{xcolor}
\usepackage{mathtools}
\usepackage{fullpage}
\usepackage{cases}
\usepackage{graphics}
\usepackage{tikz-cd}
\usepackage{caption, subcaption}
\usepackage{harpoon}
\usepackage[utf8]{inputenc}
\usepackage[top=1 in,bottom=1in, left=1 in, right=1 in]{geometry}
\usepackage{cleveref}
\usepackage{booktabs}
 \usepackage{diagbox}
\newcommand{\PP}{\mathbb{P}}
\newcommand{\RR}{\mathbb{R}}

\newcommand{\CC}{\mathbb{C}}

\renewcommand{\phi}{\varphi}

\newcommand{\rk}{\textrm{Rank}}

\DeclareMathOperator{\Supp}{Supp}

\def\frk{\frak}               

\def\Phi{{\frk n}}
\def\Phi{{\frk N}}
\def\MC{{\mathcal C}}
\def\MH{{\mathcal H}}

\def\MX{{\mathcal X}}
\def\MY{{\mathcal Y}}
\def\MZ{{\mathcal Z}}

\newcommand{\R}{\mathbb{R}}

\def\opn#1#2{\def#1{\operatorname{#2}}} 
\opn\chara{char} \opn\length{\ell} \opn\pd{pd} \opn\rk{rk}
\opn\projdim{proj\,dim} \opn\injdim{inj\,dim} \opn\rank{rank}
\opn\depth{depth} \opn\grade{grade} \opn\height{height}
\opn\embdim{emb\,dim} \opn\codim{codim}

\opn\Tr{Tr} \opn\bigrank{big\,rank}
\opn\superheight{superheight}\opn\lcm{lcm}
\opn\trdeg{tr\,deg}
\opn\reg{reg} \opn\lreg{lreg} \opn\ini{in} \opn\lpd{lpd}
\opn\size{size} \opn\sdepth{sdepth}
\opn\link{link}\opn\fdepth{fdepth}\opn\lex{lex}
\opn\div{div} \opn\Div{Div} \opn\cl{cl} \opn\Cl{Cl}
\opn\Spec{Spec} \opn\Supp{Supp} \opn\supp{supp} \opn\Sing{Sing}
\opn\Ass{Ass} \opn\Min{Min}\opn\Mon{Mon}

\newcommand{\introthmname}{}
\newtheorem{introthminn}{\introthmname}

\theoremstyle{definition}
\newtheorem{Theorem}{Theorem}[section]

\newtheorem{Definition}{Definition}[section]

\newtheorem{Example}{Example}[section]

\newtheorem{Remark}{Remark}[section]
\newtheorem{Conjecture}{Conjecture}[section]

\numberwithin{equation}{section}

\definecolor{DarkRed}{RGB}{200,20,20}
\definecolor{DarkGreen}{RGB}{0,120,0}
\definecolor{SkyBlue}{rgb}{0.16, 0.32, 0.75}
\hypersetup{
    colorlinks=true,
linkcolor = blue,    linkbordercolor=DarkRed,
     citecolor = teal,
    }

\newcommand\Perp{\protect\mathpalette{\protect\independenT}{\perp}}
\def\independenT#1#2{\mathrel{\rlap{$#1#2$}\mkern2mu{#1#2}}}
\newcommand{\ind}[3][]{\left.#2 \Perp\!_{#1}\, #3 \inD}
\newcommand{\inD}[1][\relax]{\def\argone{#1}\def\temprelax{\relax}
  \ifx\argone\temprelax\right.\else\,\middle|#1\right.{}\fi}

\title{Applications of 
Tensors in statistics and rigidity theory}
\author{Fatemeh Mohammadi}
\date{}

\begin{document}
\begin{abstract}
This is a short report on the discussions of appearance of tensors in algebraic statistics and rigidity theory, during the semester ``AGATES: Algebraic Geometry with Applications to TEnsors and Secants". We briefly survey some of the existing results in the literature and further research directions.
We first provide an overview of algebraic and geometric techniques in the study of conditional independence (CI) statistical models. 
We study different families of algebraic varieties arising in statistics. This includes the determinantal varieties related to CI statements with hidden random variables.  Such statements correspond to determinantal conditions on 
the tensor of joint probabilities of events involving the observed random variables.
We show how to compute the irreducible decompositions of the corresponding CI varieties, which leads to finding further conditional dependencies (or independencies) among the involved random variables. 
As an example, we show how these methods can be applied to extend the classical intersection axiom for CI statements. 
We then give a brief overview about secant varieties and their appearance in the study of mixture models.
We focus on examples and briefly mention the connection to rigidity theory which is further discussed in \cite{cruickshank2023identifiability}.  
\end{abstract}

\maketitle



\section{Introduction}

This note is focused on the study of mixture models and the conditional independence (CI) models with hidden variables, and their connections to tensors, determinantal and secant varieties. 



\medskip
Conditional independence provides an important tool in statistical modelling~\cite{Studeny05:Probabilistic_CI_structures}, as
it gives an interpretation to Markov fields and graphical models; see~e.g.~\cite{MDLW18:Handbook_of_graphical_models} and the references therein. These statistical models are well-studied in algebraic statistics, as they can be interpreted as solutions of polynomial equations (or inequalities). We recommend the books \cite{DrtonSturmfelsSullivant09:Algebraic_Statistics, Sullivant} for a more detailed overview of the subject. The CI models are also strongly linked to combinatorics and lattice geometry \cite{andersson1993lattice, caines2022lattice, mohammadi2018generalized}. In particular, the lattice conditional independence (LCI) models are a special type of statistical graphical models introduced in \cite{andersson1993lattice} in the context of linear Gaussian models. The idea is that in the case of a distributive lattice of linear subspaces associated with the marginal models, all CI statements could be classified via the intersections (meets) on the lattice. These models can equivalently be described via a transitive directed acyclic graph, where the conditional independences are given in terms of conditioning on ancestors in the graph \cite{caines2022lattice}. Algebraically, they are characterized as (generalized) Hibi ideals \cite{Hibi1,ene2011monomial}.

\medskip

Given a collection of random variables and knowledge of the conditional dependencies, or independencies, among them, we can ask what are the probability distributions that satisfy them. In a more general setting, some of the random variables appearing in a CI model can be prescribed as unobserved (or hidden). Then the main question is to determine when certain constraints on the observed variables arise from conditions on the hidden variables \cite{Steudel-Ay}. This problem can be restated algebraically by noting that probability distributions satisfying CI statements are the solutions of certain polynomial equations \cite{DrtonSturmfelsSullivant09:Algebraic_Statistics, Sullivant} which generate the so-called CI ideal. The distributions satisfying a given collection of CI statements can be recovered by intersecting the CI ideal with the probability simplex. When there are no hidden variables, these polynomials are binomials and their associated ideals are well-studied; see e.g.~\cite{Fink,herzog2010binomial,Rauh,SwansonTaylor11:Minimial_Primes_of_CI_Ideals}. However, in the presence of hidden variables, the polynomials become far more complicated of arbitrarily high degrees and very difficult to calculate; see e.g.~\cite{pfister2021primary,clarke2020conditional, MatroidsCIStatements}.

\begin{Example}
Let $X, Y, Z$ be discrete random variables taking finitely many values in the sets $\MX, \MY$ and $\MZ$, respectively. The joint probability distribution of $X, Y$ and $Z$ can be identified with a $3$-dimensional  
tensor 
\[
P=(p_{x,y,z})_{x\in\MX, y\in\MY, z\in\MZ},\quad\text{where}\quad p_{x,y,z} = \PP(X = x, Y = y, Z = z).
\]
Each CI statement is equivalent to bounding the rank of slices of this tensor. 
More precisely, the variables $X$ and $Y$ are \textit{independent} given $Z$, denoted $\ind XY \mid Z$, if and only if for each $z \in \MZ$ the matrix $P_z:=(p_{x,y,z})_{x\in\MX, y\in\MY}$ has rank one. 
Suppose that $\ind XY \mid Z$ where $Z$ is a hidden variable. So the observed distribution is the joint distribution of $X$ and $Y$ which is the marginal tensor $P^{X,Y} = \sum_{z \in \MZ} P_z$. This is called a \textit{flattening} of the tensor $P$. By the CI statement $\ind XY \mid Z$, we have that each $z$-slice $P_z$ is a rank-one matrix. The joint distribution of $X,Y$ is a sum of these matrices, hence $P^{X,Y}$ has rank at most $|\MZ|$, or equivalently, all $(|\MZ|+1)$-minors of $P^{X,Y}$ vanish. 
\end{Example}

This motivates the definition of the \textit{CI ideal} for models with  hidden random variables, that we  study in Section~\ref{sec:CI_ideals}. For a concrete instance of the above example, see Example~\ref{example:CI-3-4-2}.

\medskip
In Section~\ref{sec:matroids}, we give an overview of our results from \cite{MatroidsCIStatements}, where motivated by the realizability and causality problems in statistics, we study the associated varieties of hypergraphs from the projective geometry and matroid theory viewpoints. In particular, we describe the connections between  
conditional independence (CI) models 
in statistics \cite{Studeny05:Probabilistic_CI_structures,DrtonSturmfelsSullivant09:Algebraic_Statistics,Sullivant}, projective geometry \cite{richter2011perspectives,lee2013mnev}, and the theory of matroids \cite{Oxley, piff1970vector} and their realization spaces 
\cite{mnev1985manifolds, mnev1988universality, sturmfels1989matroid}. 
We will focus on the specific family of grid hypergraphs, whose 
corresponding varieties have explicit interpretations in terms of matroids.  
We will show how the hypergraph ideals are related to CI ideals. Our goal is to compute primary decompositions of these ideals; see Example~\ref{exa:k2l4s2t3} and Theorem~\ref{MatroidsCIStatements}. 

\medskip
Finally, we provide a brief overview about secant varieties and their connections to mixture models.
We focus on examples and refer to \cite{DrtonSturmfelsSullivant09:Algebraic_Statistics} for more details. In particular, the study of mixture models
is related to non-negative matrices of low-rank. They are also connected to local and global rigidity of bar-joint frameworks; see~e.g.~\cite{krone2021uniqueness} and \cite{cruickshank2023identifiability}.  
\section{Conditional independence ideals and the intersection axiom}\label{sec:CI_ideals} 
Given a collection of observed random variables $X_1, \dots$, $X_n$ and a hidden random variable $H$ taking values in $\MX_1, \dots, \MX_n, \MH$, respectively, we define the polynomial ring $R = \CC[p_{i_1, \dots, i_n} : i_1 \in \MX_1, \dots, i_n \in \MX_n]$ which has one variable for each outcome in the joint probability distribution of the observed variables. Let $\MC$ be a collection of CI statements for the random variables $X_1, \dots$, $X_n$ and $H$. The \textit{CI ideal} $J_\MC \subseteq R$ is defined to be the ideal generated by rank constraints on the (flattenings of the) joint probability tensor $P=(p_{i_1, \dots, i_n})$.

\begin{Example}\label{example:CI-3-4-2}
{\rm 
Consider the observed random variables $X, Y_1, Y_2$ and hidden random variables $H_1, H_2$. 
Consider the following collection of CI statements:
\begin{eqnarray}\label{eq:C}
\MC : \ \ 
\ind X {Y_1} \mid \{Y_2, H_1\} \quad \textrm{and} \quad  \ind X {Y_2} \mid \{Y_1, H_2\}.
\end{eqnarray}
Let 
$\MX= \MY_1 = \{1,2,3\}, \MY_2 = \{1,2,3,4\}$ and $\MH_1 = \MH_2 = \{0,1\}$.  The joint distribution of $Y_1$ and $Y_2$ has state space $\MY = \MY_1 \times \MY_2$ which is identified with the $3 \times 4$ matrix $\MY$ with values in the set $[12]$.
Since the hidden variables $H_1$ and $H_2$ take two different values, the CI ideal $J_{\MC}$ is generated by $3$-minors of the matrix of variables $P = (p_{x,y})_{x \in \MX, y \in \MY}$. Explicitly, the CI ideal $J_{\ind{X}{Y_1}\mid\{Y_2, H_1 \}}$ is generated by the four $3$-minors of $P$ whose columns are indexed by $C_1 = \{1,2,3\}, \dots, C_4 = \{10,11,12 \}$, i.e.~the columns of $\MY$. Similarly, $J_{\ind{X}{Y_2}\mid\{Y_1, H_2 \}}$ is generated by the $3$-minors of $P$ whose columns are indexed by $3$-subsets of the rows $R_1, R_2$ and $R_3$ of $\MY$. In Figure~\ref{fig:flatten_tensor}, we show the flattening of the probability tensor for each $Y_2$-slice that gives rise to the $3$-minors that generate the CI ideal.
In Section~\ref{sec:matroids}, we identify the CI ideals $J_{\mathcal{C}}$ as a hypergraph ideal; see Example~\ref{ex:andreas}. 
In particular, we will show that the ideal $J_{\mathcal{C}}$ has two prime components that can be identified as the ideals associated to certain matroids.} 
Note that $H_1$ and $H_2$ are hidden random variables taking the same number of values, and so the ideals do not distinguish them.

\begin{figure}[h]
    \centering
    \includegraphics[scale=0.7]{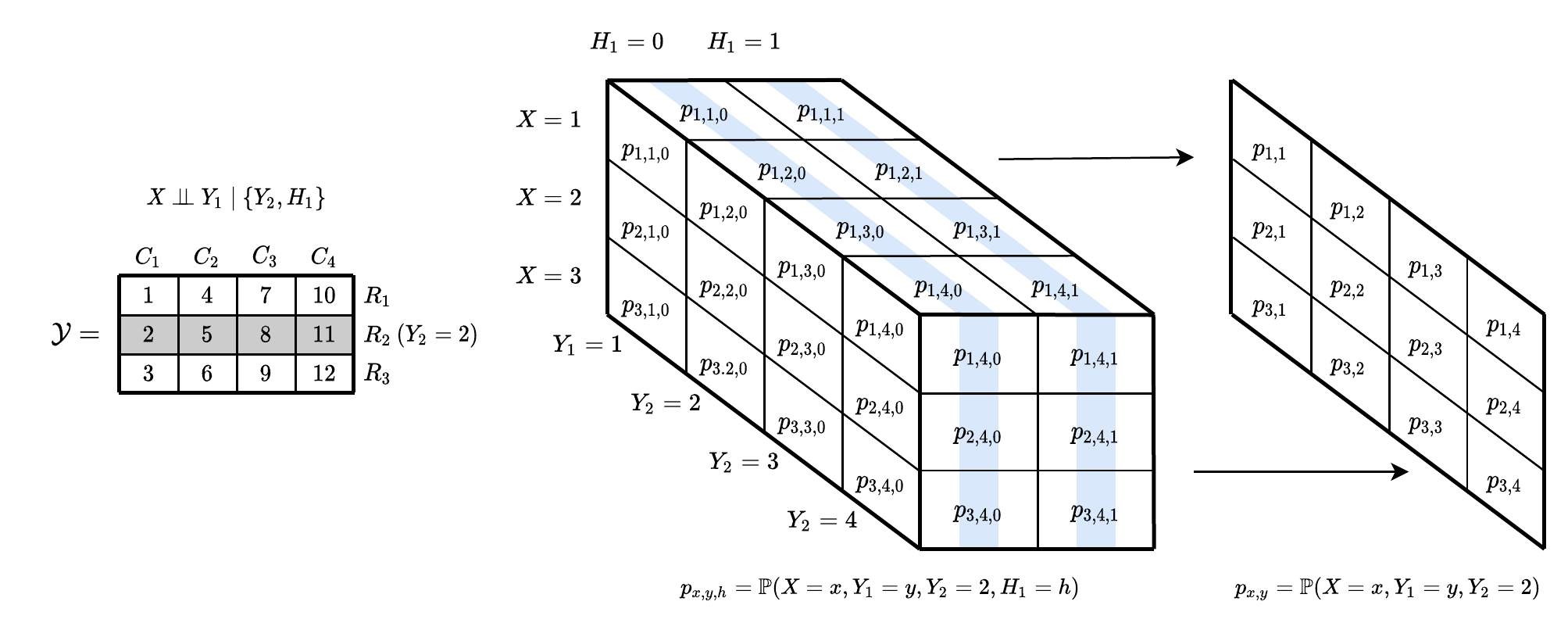}
    \caption{Depiction of the flattening of the joint probability tensor $P$ in Example~\ref{example:CI-3-4-2} that gives rise to the CI ideal $J_{\ind{X}{Y_1} \mid \{Y_2, H_1 \}}$. Consider the slice $Y_2 = 2$ of $P$. The $H_1 = 0$ and $H_1 = 1$ slices of the tensor $(p_{x,y,h})$, shown in the figure with stripes, each have rank one. Therefore the resulting flattening given by $p_{x,y} = p_{x,y,0} + p_{x,y,1}$ has rank at most two. Therefore the $3$-minors of the matrix $(p_{x,y})$ are zero.
    }
    \label{fig:flatten_tensor}
\end{figure}
\end{Example}

\noindent{\bf Intersection axiom.} Suppose that $H_1$ and $H_2$ are both constant. The \textit{intersection axiom} states that any distribution satisfying $\MC$ in \eqref{eq:C} \textit{generically} satisfies $\ind X {\{Y_1, Y_2\}}$, that is the distributions have non-zero probabilities. Consider the CI ideal $J_\MC$, which is the intersection of prime ideals \cite{Fink}. One of these prime components is given by $J_{\ind{X}{\{Y_1, Y_2\}}}$. In terms of varieties, this is the component containing all distributions that are \textit{fully supported}, i.e.~do not have any \textit{structural zeros}. The other components contain structural zeros which can be combinatorially classified \cite{herzog2010binomial, Rauh}.

\begin{Remark}
The intersection axiom has been studied for various CI models, see e.g.~\cite{herzog2010binomial, Rauh, clarke2021conditional, pfister2021primary}. In fact all the previous cases can be identified as ideals of (grid) hypergraphs; see \cite{MatroidsCIStatements}. In particular, in each case the corresponding ideal 
has a distinguished prime component with a particular statistical significance, since the distributions which lie inside do not contain any structural zeros \cite{Sturmfels02:Solving_polynomial_equations}. 
This prime ideal can be realized as the CI ideal of $\ind X {\{Y_1, Y_2\}} \mid H_2$. We may therefore deduce a hidden variable version of the intersection axiom as follows:
\[
\MC = \{\ind X {Y_1} \mid \{H_1,Y_2\}, \ \ind X {Y_2} \mid \{Y_1, H_2\} \}
\implies
\ind{X}{\{Y_1, Y_2\} \mid H_2}.
\]
\end{Remark}


We end this section with a related conjecture by Matúš \cite{matu1999conditional}, and explain the connection to realization spaces of matroids in the following section.

\begin{Conjecture}[\cite{matu1999conditional}]\label{conj:matus_rational}
For any discrete conditional independence model $\MC$, there exists a distribution $p$ in the zero set of $J_{\MC}$ such that all joint probabilities of $p$ are rational. 
\end{Conjecture}

\section{Hypergraph varieties and their matroid stratifications}\label{sec:matroids}

Let $\mathbb{K}$ be a field, $d\leq n$ be two positive integers, 
$X = (x_{ij})$ be a $d\times n$ matrix of indeterminates and $R =\mathbb{K}[X]$ be the polynomial ring over $\mathbb{K}$ in the indeterminates $x_{ij}$. We write determinants of submatrices of $X$ as $[I | J]_X$ where $I$ and $J$ are respectively the sets of rows and columns of the submatrix. 
We denote by $x_i$ the $i^{\rm th}$ column of $X$ and by $X_F$ the submatrix of $X$ with columns indexed by $F\subseteq[n]$. We recall the following definition from \cite{Fatemeh, Fatemeh2}.

\begin{Definition}
\label{def:hypergraph}
A (simple) hypergraph $\Delta$ on the vertex set
$[n]$ is a subset of the power set $2^{[n]}$. We assume that no proper subset of an element of $\Delta$ is in $\Delta$.
The elements of $\Delta$ are called (hyper)edges.  
\begin{itemize}
    \item The \emph{determinantal hypergraph ideal} of $\Delta$ is
  \begin{equation*}
    I_{\Delta}= \big\langle [A|B]_X : A\subseteq [d], B\in\Delta, |A| = |B|  \big\rangle \subset R.
  \end{equation*}
  \item The variety of $\Delta$ is the zero set of $I_{\Delta}$ which is given by 
  \[
  V_{\Delta}=\{X\in \mathbb{C}^{d\times n}:\ \rk(X_F) < |F| \text{ for each }F \text{ in } \Delta\}.
  \] 
  \end{itemize}
\end{Definition}

In particular, we have that:

\begin{Theorem}[\cite{MatroidsCIStatements}]\label{MatroidsCIStatements} 
Every hypergraph variety $V_\Delta$ is the union of matroid varieties. The union is taken over all realizable matroids whose dependent sets contain the edges of $\Delta$.
\end{Theorem}
Note that many of these components might be redundant or reducible. Finding a minimal irreducible decompositions of hypergraph varieties is a challenging open problem.

\medskip
We now recall the definitions of the realization space of a matroid and its associated variety. We refer the reader to \cite{Oxley} for basic definitions concerning matroids.
\begin{Definition}\label{def:prelim_realisation}
Let $M$ be a matroid on $[n]$ of rank $r$ and let $d \ge r$. If $\mathbb{K}$ is a field, a \emph{realization} of $M$ in $\mathbb{K}^d$ is a collection of vectors $X = \{x_1, \dots, x_n \} \subset \mathbb{K}^d$ such that 
\[
\{x_{i_1}, \dots, x_{i_p}\} \subset X \text{ is linearly dependent} \iff \{i_1, \dots, i_p \} \text{ is a dependent set of } M.
\]
If such a collection of vectors exists, we say that the matroid is \emph{realizable} over $\mathbb{K}$. 
Here, we focus on the case $\mathbb{K}=\CC$.
The \emph{realization space} of $M$ in $\CC^d$ is
\[
\Gamma_{M} = \{X \subset \CC^d : X \text{ is a realization of } M \}.
\]
\end{Definition}
Each element of $\Gamma_{M}$ can be identified with a $d \times n$ matrix $X$ over $\CC$.
The \emph{matroid variety} $V_{M} = \overline{\Gamma_{M}}$ is the Zariski closure of the realization space of $M$, and $I_{M} = I(V_{M}) \subseteq \CC[X]$ denotes its corresponding ideal.


\begin{figure}[h]
    \centering
    \includegraphics[scale = 0.8]{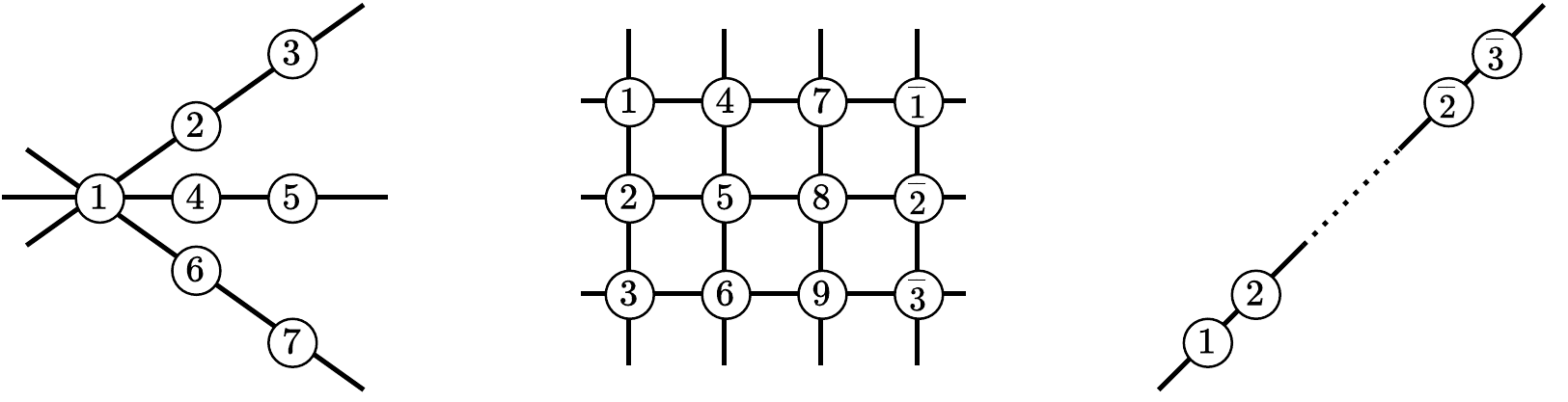}
    \caption{(Left) A realization of $M$ from Example~\ref{ex:combthreelines}.
    (Center and Right) Realizations of matroids corresponding to the prime components of $I_{\Delta}$ in Example~\ref{ex:andreas}.
    }
    \label{threelines}
\end{figure}

\begin{Example}\label{ex:combthreelines}

Let $d = 3$, $n = 7$ and $\Delta = \min(\{123, 145, 167 \} \cup \binom{[7]}{4})$. It is easy to check that $\Delta$ is the collection of circuits for a matroid $M$, hence  $M$ is the unique minimal matroid for $\Delta$.  
The associated ideal 
is the ideal $I_{\Delta} = \langle [123], [145], [167] \rangle \subseteq \CC[x_{1,1}, \dots, x_{3,7}]$.
with two prime components:
\[
I_{\Delta} = I_1 \cap I_2 =
\langle x_{1,1}, x_{2,1}, x_{3,1} \rangle \cap 
\langle [123], [145], [167], [234][567] - [235][467] \rangle 
\]
The ideal $I_1$ is the ideal of the matroid $M'$ that has a single circuit $1$ and has obtained from $M$ by setting $1$ to be a loop. The ideal $I_2$ is the ideal of the matroid $M$ and the generator $[234][567] - [235][467]$ of $I_2$ is a geometric condition satisfied by six generic points lying on three lines that intersect at a common point, as shown in Figure~\ref{threelines} (left).



\end{Example}

\begin{Example}
\label{ex:andreas}
Let $d = 3$ and $E=[9]\cup\{\bar{1},\bar{2},\bar{3}\}$. Let  $\Delta$ be the hypergraph
$$\Delta=\{123, 456, 789,\bar{1}\bar{2}\bar{3},
147,\bar{1}14,\bar{1}47,\bar{1}17, 258,\bar{2}58,\bar{2}28,\bar{2}25,369,\bar{3}69,\bar{3}39,\bar{3}36\}$$ depicted in Figure~\ref{threelines} (center). It is shown in \cite[Theorem~4.1]{pfister2021primary} that $I_{\Delta}$ has two prime components:
\begin{itemize}
\item The first component is generated by all $3$-minors of a generic $3\times 12$ matrix. Hence, it corresponds to the matroid $M_0$ with dependent sets, $$\mathcal{C}(M_0)=\{\text{all $3$-subsets of $E$}\}.$$

\item The second component is generated by $44$ polynomials of which $16$ are the original generators of $I_{\Delta}$
and the remaining $28$ generators are all homogeneous of degree $12$, which are obtained by the geometric constraints of \textit{quadrilateral sets}; see \cite[\S8]{richter2011perspectives}.
\end{itemize}
\end{Example}

We now generalize the above example to so-called grid hypergraphs $\Delta^{s,t}$. 
\medskip

\noindent{\bf Grid hypergraphs $\Delta^{s,t}$.}
Consider the $k \times \ell$ matrix of integers
\begin{eqnarray}\label{eq:Y}
\MY = 
\begin{bmatrix}
    1       & k + 1     & \dots     & (\ell - 1)k + 1 \\ 
    2       & k + 2     & \dots     & (\ell - 1)k + 2 \\
    \vdots  & \vdots    & \ddots    & \vdots \\
    k       & 2k        & \dots     & \ell k
\end{bmatrix} \ .
\end{eqnarray}
Then for each $1 \le i \le k$, the rows of $\MY$ are denoted
\[
R_i = \{ \MY_{i,1}, \MY_{i,2}, \dots, \MY_{i, \ell} \} = \{i, k+i, \dots, (\ell - 1)k + i \}
\] and for each $1 \le i \le \ell$, the columns of $\MY$ are denoted
\[
C_j = \{ \MY_{1,j}, \MY_{2,j}, \dots, \MY_{k,j}\} = \{(j-1)k + 1, (j-1)k + 2, \dots, (j-1)k + k \}.
\]
For each $s$ and $t$ (with $s \le k$, $t \le \ell$), we define $\Delta^{s,t}$ to be the following collection of subsets of $[k\ell]$, 
\[
\Delta^{s,t} = \bigcup_{1 \le i \le k} \binom{R_i}{t} \cup \bigcup_{1 \le j \le \ell} \binom{C_j}{s}.
\]

\begin{figure}[h]
    \centering
    \includegraphics[scale=0.25]{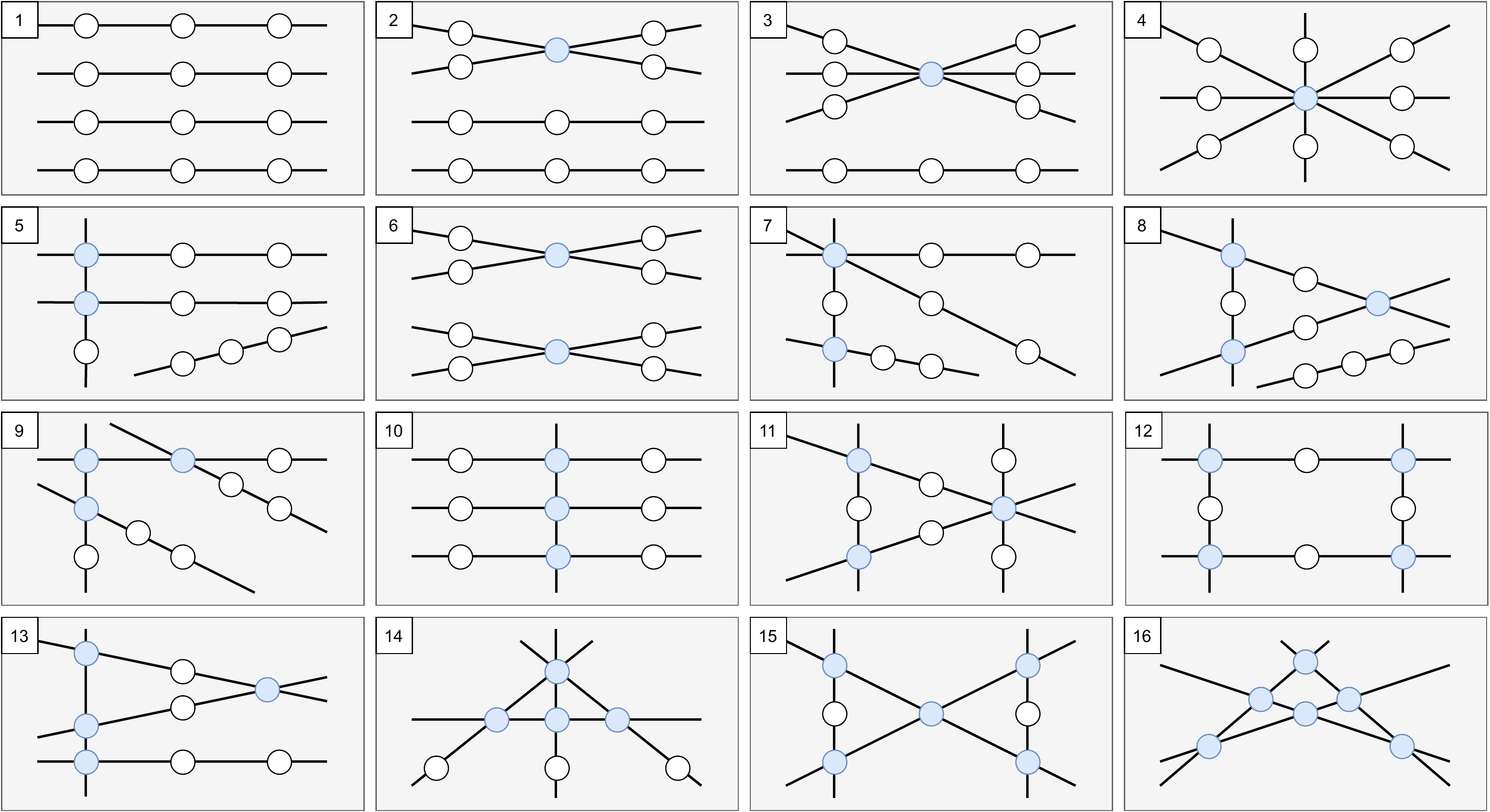}
    \caption{The combinatorial types of point and line arrangements with $4$ lines.}
    \label{fig:4xn_combinatorial_types}
\end{figure}

\begin{Example}
\label{exa:k2l4s2t3}
Let $k = 4, \ell = 7, s = 2$ and $t = 3$. We have
\[
\MY = 
\begin{bmatrix}
1 & 5 & 9  & 13 & 17 & 21 & 25 \\
2 & 6 & 10 & 14 & 18 & 22 & 26 \\
3 & 7 & 11 & 15 & 19 & 23 & 27 \\
4 & 8 & 12 & 16 & 20 & 24 & 28
\end{bmatrix},
\quad
\Delta^{2,3} = \left\{ 
\binom{\{1,5,9,13,17,21,25\}}{3} 
\cup
\dots 
\cup
\binom{\{25,26,27,28 \}}{2}
\right\}.
\]
Calculating the dependent matroids for $\Delta^{2,3}$, we find that all such matroids are matroids of 
line arrangements. There are $10$ combinatorial types of line arrangements which appear as these matroids. Those with exactly $4$ distinct lines are drawn in Figure~\ref{fig:4xn_combinatorial_types}, specifically these are types: 14, 15, 16. The $7$ remaining combinatorial types are point and line arrangements with at most $3$ lines.
\end{Example}

Grid hypergraphs provide important families of examples, as for every grid hypergraph there is a unique minimal matroid containing it, and conversely, for every matroid $M$ there exists a grid hypergraph and a dependent matroid whose restrictions is isomorphic to $M$.
\begin{Theorem}[\cite{MatroidsCIStatements}]\label{thm:s-t-3}
Let $s,t,k,\ell,d$ be positive integers such that $3\leq s\leq t\leq\ell$, $s\leq k$, and $t\leq d\leq s+t-3$. Then $\MC=\min(\Delta^{s,t}\cup\binom{[k\ell]}{d+1})$ is the collection of circuits of an $\mathbb{R}$-realizable matroid on $[k\ell]$ of rank $d$. This is the unique minimal matroid for $\Delta^{s,t}$ in this case.
On the other hand, for every matroid $M$, there exists a grid hypergraph $\Delta^{s,t}$ and a dependent matroid $M'$ for $\Delta^{s,t}$ such that a \textit{restriction} of $M'$ is isomorphic to $M$.
\end{Theorem}

\noindent{\bf Correspondence of $\Delta^{s,t}$ with CI models.}
The minimal prime decomposition of the ideal of $\Delta^{s,t}$ has been extensively studied in \cite{herzog2010binomial, Rauh, clarke2020conditional,pfister2021primary}. 
In each of these cases, there is a straightforward matroidal description of the prime components. The hypergraph ideals of $\Delta^{s,t}$ arise in the following setting:
Consider three observed variables $X,Y_{1},Y_{2}$, taking values in the finite sets $\MX,\MY_{1},\MY_{2}$ of cardinalities $|\MX|=d$, $|\MY_{1}|=k$, $|\MY_{2}|=\ell$, and two hidden variables $H_{1},H_{2}$, taking values in the finite sets $\MH_{1},\MH_{2}$ of cardinalities $|\MH_{1}|=s-1$,  $|\MH_{1}|=t-1$.
The joint distribution of the observed variables can be identified with a non-negative matrix
$P\in\R^{\MX\times\MY}$, where $\MY=\MY_{1}\times\MY_{2}$.  The matrix $P$ has $d$ rows and
$k\ell$ columns, and its entries $p_{x,(y_{1},y_{2})}$ sum to one.  The row indices
correspond to states $x\in\MX$, and the column indices correspond to joint states
$(y_1,y_{2})\in\MY$.

\medskip

We now explain the potential approach to study Conjecture~\ref{conj:matus_rational}. In Theorem~\ref{thm:s-t-3}, we have seen that for large enough $s,t,k,\ell$, any matroid may appear among the dependent matroids for $\Delta^{s,t}$. A natural approach is to carefully choose additional conditional independence and dependence statements for the model \eqref{eq:C}, in order to guarantee that any distribution $p$ satisfying $\MC$ is a realization of a given, realizable, matroid. 
Note that there exist matroids that are not realizable over the rationals but are realizable over a real field extension. Hence, this might lead to a characterization of CI models with hidden variables for which Conjecture~\ref{conj:matus_rational} does not hold.

\medskip

\begin{Example}
Let $k=4,\ell=3,s = 2, t = 3$. 
The prime components of the ideal of $I_{\Delta^{s,t}}$ can be identified as the ideals associated to configurations with 12 points and $4$ lines. With $12$ points, we observe all combinatorial types of line arrangements with at most $4$ lines, depicted in Figure~\ref{fig:4xn_combinatorial_types}. The ideal of such configuration shows up as a prime component of $I_{\Delta^{s,t}}$.  
\end{Example}

\noindent{\bf Generic Rigidity.}\label{sec:applications_rigidity}
We now briefly outline the well-studied application of matroid theory in rigidity theory and its connection to CI models.  For a comprehensive text on the topic see \cite{sitharam2017handbook}. 

 \smallskip
Rigidity theory of bar and joint
frameworks investigates the sufficient and necessary conditions for the uniqueness of a point configuration given some of
the pairwise distances. Let $K_n$ denote the complete graph on $n$-vertices with the edge set $E(K_n)$. We associate to it the \textit{generic $d$-rigidity matroid} $R_{d,n}$ with ground set $E(K_n)$. The matroid $R_{d,n}$ is defined by a realization of a generic embedding of $K_n$ into $\RR^d$, given by the \textit{rigidity matrix}. The matroid $R_{d,n}$ holds important information about the generic rigidity of bar and joint frameworks, and characterizing the matroid $R_{d,n}$ is an important open problem in rigidity theory.

 \smallskip
 
Recently, Jackson and Tanigawa \cite{jackson2021maximal} have considered \textit{abstract $d$-rigidity matroids}, originally studied by Graver \cite{graver91rigidity}, which can be defined as the family of matroids $M$ on ground set $E(K_n)$ such that all $K_{d+2}$ subsets are circuits and $M$ has rank $dn - \binom{d+1}{2}$. It is an open problem to determine whether this family has a unique minimally dependent matroid. In particular, if $d = 3$ then it is open as to whether this matroid is $R_{d,n}$. We hope to apply the results of this session to explore algorithmic approaches to this problem.
In particular, for the grid hypergraphs $\Delta^{s,t}$, when $d \le s + t - 3$, we have seen that there is a unique minimal matroid $M^{s,t}$ for $\Delta^{s,t}$. Analogously to rigidity matroids, we can think of $M^{s,t}$ as the unique minimally dependent matroid with circuits given by $\Delta^{s,t}$ and with rank $d - 1$. 

\smallskip
Another related avenue of research is the study of algebraic matroids and matrix completion. 
\begin{Definition}
Let $V \subseteq \CC^n$ be an irreducible affine variety. For each subset $S \subseteq [n]$, we define the projection map $\pi_S : \CC^n \rightarrow \CC^S$ which takes each point $(x_i)_{i \in [n]}$ to the point $(x_i)_{i \in S}$. The dependent sets of the algebraic matroid of $V$ are defined to be the subsets $S \subseteq [n]$ such that Zariski closure $\overline{\pi_S(V)} \neq \CC^S$ is not the entire space.
\end{Definition}

We note that the matroid $M^{s,t}$ arises as an \textit{algebraic matroid} in the following sense. 

\begin{Example} Fix $s,t,k,\ell,d$ positive integers and $Y = (y_{i,j})$ a $k\times\ell$ matrix of variables.
Let $V$ be the affine variety whose points are all $k \times \ell$ matrices of rank at most $d$, such that each column contains at most $s$ nonzero entries and each row contains at most $t$ nonzero entries. Explicitly, the variety $V$ is the vanishing set of the following ideal in the polynomial ring $\CC[Y] = \CC[y_{i,j} : i \in [k], j \in [\ell]]$:
\[
\left\langle [I|J]_Y : I \in \binom{[k]}{d}, J \in \binom{[\ell]}{d} \right\rangle
+
\left\langle 
\Pi_{i \in I} y_{i,j} : I \in \binom{[k]}{s}, j \in [\ell]
\right\rangle
+
\left\langle
\Pi_{j \in J} y_{i,j} : i \in [k], J \in \binom{[\ell]}{t}
\right\rangle.
\]
It is not immediately obvious that $V$ is an irreducible variety. However, if we assume $V$ is irreducible, then it is easy to see that each edge of $\Delta^{s,t}$ is a dependent set of the algebraic matroid of $V$. Furthermore, if $d \le s + t -3$ then, by Theorem~\ref{thm:s-t-3}, we have that the algebraic matroid of $V$ coincides with the matroid $M^{s,t}$ above.
\end{Example}



We note that the use of hypergraphs to characterize rigid structures is analogous to understanding irreducible structures in our setting. 
In \textit{scene analysis} \cite{Whiteley1989hypergraph_scene_analysis}, given a \textit{plane picture} which contains the data of a projection of given incidence structure, the central question is to determine \textit{generic} preimages of this projection, i.e.~find generic realizations of an incident structure subject to a condition on its projection. In our setting, this translates to finding the irreducible components of the (Zariski closures of the) realizations of these incidence structures.



\section{Secant varieties and mixture models}\label{sec:secant}
Secant varieties are classical objects in algebraic geometry. They are defined using the notion of join of two varieties.
Let $\mathbb{K}$ be a field, and let $U$ and $V$ be two affine varieties. Then the join of $U$ and $V$ is the following variety:
\[
\mathcal{J}(U,V)=\overline{\{\lambda u+(1-\lambda)v:\ u\in U,\ v\in V,\ \lambda\in \mathbb{K}\}},
\]
which is the Zariski closure of the set of all points on the lines connecting a point in $U$ and a point in $V$. The secant varieties of $V$ are defined inductively as follows:
\[
{\rm Sec}^2(V):=\mathcal{J}(V,V)\quad{\rm and}\quad {\rm Sec}^k(V):=\mathcal{J}({\rm Sec}^{k-1}(V),V).
\]
Given two subsets $U$ and $V$ of $\mathbb{R}^n$ we can define their mixture as:
\[
{\rm Mixt}(U,V)=\{\lambda u+(1-\lambda)v:\ u\in U,\ v\in V,\ 0\leq\lambda\leq 1\},
\]
which is the set of all convex combinations of a point in $U$ and a point in $V$. We also define:
\[
{\rm Mixt}^2(V):={\rm Mixt}(V,V)\quad{\rm and}\quad {\rm Mixt}^k(V):={\rm Mix}({\rm Mixt}^{k-1}(V),V).
\]
For a semi-algebraic set $V$, the secant variety of its Zariski closure coincides with the Zariski closure of its mixture, i.e.~${\rm Sec}^k(\overline{V})=\overline{{\rm Mixt}^k(V)}$. However, ${\rm Sec}^k(\overline{V})$ and ${\rm Mixt}^k(V)$ can be in general very different; see~e.g. Example~4.1.2 in \cite{DrtonSturmfelsSullivant09:Algebraic_Statistics}.

\medskip

Mixture models are used to construct complex statistical models from the simple ones. 
Here, we only provide an example and we recommend \cite{DrtonSturmfelsSullivant09:Algebraic_Statistics} for a more detailed overview of the subject.
\begin{Example}
Let $X$ and $Y$ be two independent random variables taking finitely many values in the sets $\mathcal{X}$ and $\mathcal{Y}$, respectively. The independence model $\mathcal{M}_{\ind X {Y}}$ 
is the set of all non-negative rank $1$ probability matrices of size $|\mathcal{X}|\times|\mathcal{Y}|$, and the $k^{\rm th}$ mixture model of $\mathcal{M}_{\ind X {Y}}$, denoted by Mixt$^k(\mathcal{M}_{\ind X {Y}})$ is the set of probability matrices of non-negative rank at most $k$.
\end{Example}
Some of the main algebraic problems about mixture models with ample applications in statistical inference are the following: (1) computing the dimensions of these models (2) determining their singularities (3) and characterizing identifiable models among them. The latter means given a probability distribution in the model, we want to determine whether
the parameters of the model can be reconstructed. The dimensions of these models have been also studied in comparison with the dimension of the secant varieties. 

\smallskip
\noindent{\bf Connections with rigidity theory.}
Rigidity theory provides useful techniques in the study of the uniqueness of nonnegative matrix factorizations \cite{krone2021uniqueness} and the uniqueness of low-rank matrix completion \cite{singer2010uniqueness}. See  \cite{cruickshank2023identifiability} for the connections of rigidity theory, matrix completions, and the expected dimensions of secant varieties.

\medskip
\noindent{\bf Acknowledgement.}
The author would like to thank Jaroslaw Buczy\'nski for valuable comments on the earlier draft of this note, and Oliver Clarke and Kevin Grace for many helpful discussions on matroid varieties, and James Cruickshank, Oleg Karpenkov, Anthony Nixon and Shin-ichi Tanigawa for helpful discussions on rigidity theory during the projects \cite{ cruickshank2022global, mohammadi2022rational}.
This work is partially supported by the Thematic Research Programme ``Tensors: geometry, complexity and quantum entanglement", University of Warsaw, Excellence Initiative – Research University and the Simons Foundation Award No. 663281 granted to the Institute of Mathematics of the Polish Academy of Sciences for the years 2021-2023.  The author would
like to thank Jaroslaw Buczy\'nski, Weronika Buczy\'nska, Francesco Galuppi, and Joachim
Jelisiejew for organizing this programme.

\medskip\bibliographystyle{alpha}
\bibliography{Det.bib}

\bigskip
\noindent 
\footnotesize{\textbf{Author's addresses:}

\medskip
\noindent
Department of Computer Science, KU Leuven, Celestijnenlaan 200A, B-3001 Leuven, Belgium\\ 
   Department of Mathematics, KU Leuven, Celestijnenlaan 200B, B-3001 Leuven, Belgium 
 \\ E-mail address: {\tt fatemeh.mohammadi@kuleuven.be}
}\end{document}